%% file: Arxiv-BI.tex
\documentclass[11pt]{amsart}
\usepackage{graphicx}
\usepackage{enumerate}
\usepackage[usenames,dvipsnames,svgnames,table]{xcolor}

\theoremstyle{plain}
 \newtheorem{thm}{Theorem}[section]
      \newtheorem{cor}[thm]{Corollary}
      
      \newtheorem{lem}[thm]{Lemma}

      \newtheorem{prop}[thm]{Proposition}

      \newtheorem{rem}[thm]{Remark}

\usepackage{amstext}

\begin{document}

\title{Coherent band pathways between knots and links}

\author{Dorothy Buck}
\address{Dept of Mathematics, Imperial College London, South Kensington, London England SW7 2AZ}
\email{d.buck@imperial.ac.uk}

\author{Kai Ishihara}
\address{Yamaguchi University, 1677-1 Yoshida, Yamaguchi-shi, Yamaguchi, Japan, 753-8511}
\email{kisihara@yamaguchi-u.ac.jp}

\begin{abstract}
We categorise coherent band (aka nullification) pathways between knots and 2-component links.  
Additionally, we characterise the minimal coherent band pathways (with intermediates) between any two knots or 2-component links with small crossing number.  
We demonstrate these band surgeries for knots and links with small crossing number.  
We apply these results to place lower bounds 
on the minimum number of recombinant events separating DNA configurations, restrict the recombination pathways and determine
chirality and/or orientation of the resulting recombinant DNA molecules.
\end{abstract}

\maketitle

\section{Introduction}

\begin{figure}[htbp]
\begin{center}
\includegraphics[scale=0.8]{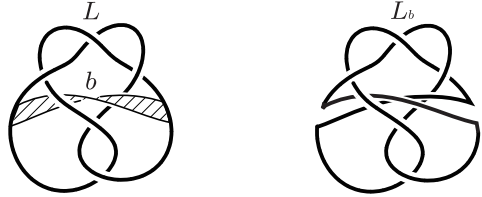}
\caption{A band surgery.}\label{Fig;band}
\end{center}
\end{figure}

Let $L$ be a link and $b$ a disk which intersects with $L$ at a disjoint union $\alpha$ of two arcs in the boundary, i.e. $L\cap b=\alpha\subset\partial b$. 
Let $\beta$ be a disjoint union of arcs obtained from $\partial b$ by removing $\alpha$, i.e. $\beta=\overline{\partial b-\alpha}$.
Then we obtain a link $L_b=(L-\alpha)\cup\beta$ by replacing $\alpha$ in $L$ with $\beta$. 
We call this operation constructing $L_b$ from $L$ a {\em band surgery} and $b$ a {\em band} of the band surgery, see Figure \ref{Fig;band}.  If $L$ and $L_b$ are oriented links,  and a band surgery $L\to L_b$ preserves the orientations of $L$ and $L_b$ (except for the band $b$), the band surgery is said to be {\em coherent}.  (Band surgeries for unoriented links are also called $H(2)$-moves \cite{Kn-Mi}.)

A coherent band surgery always changes the number of link components by one.
Conversely, a band surgery which changes the number of link components becomes a coherent band surgery by taking appropriate orientations of both links. 

Band surgeries are an important area of knot theory, and have been well-studied.  In addition to considering band surgeries in terms of cobordisms or nullification (both discussed below) Kanenobu \cite{Kn1.5}, and joint with Abe \cite{A-Kn} and Miyazawa \cite{Kn-Mi}, has
considered surgeries which preserve the number of link components as well as band surgeries on unoriented links.
Also the second author and Shimokawa gave a table of pairs of (oriented) knots and $(2,2k)$-torus links with and without coherent band surgeries \cite{I-S}.

\subsection{Cobordisms.}
Cobordisms are essential tools in examining surfaces embedded in 4-space, and are natural analogues as band surgeries.
Given two oriented links $L_0$ and $L_1$, a \textit{cobordism} from $L_0$ to $L_1$ is a pair $(S^3 \times I, F)$, where $F$ is a properly embedded oriented surface in $S^3 \times I$ such that $\partial F \cap (S^3 \times \{i\})=(-1)^{i+1}L_i$ for $i=0,1.$ 
From the cobordism perspective, a band surgery is a cobordism obtained by attaching a single 1-handle to
the surface $L_0 \times [0,\frac{1}{2}]$ in $S^3 \times I.$

If a knot has $m$ band surgeries producing an $(m+1)$-component trivial link, then the knot is called a \textit{ribbon knot}.  (Note that a ribbon knot is a slice knot, namely bounds a disc in the $4$-ball.)
The least number of such $m$ is called the \textit{ribbon fusion number}.
Kanenobu gave several conditions of coherent band surgery for two oriented links via ribbon fusion number \cite{Kn1}.

A 2-component link is said to be \textit{band-trivialisable} if there exists a band surgery producing a trivial knot. 
If a two component link is band-trivialisable, the link has $4$-ball genus zero. 
Kanenobu also applied his methods to band-trivialisability, and gave examples of $2$-component links which are $4$-ball genus zero and not band-trivialisable \cite{Kn2}.

\subsection{Nullification.}
\begin{figure}[htbp]
  \begin{center}
    \includegraphics[scale=0.2]{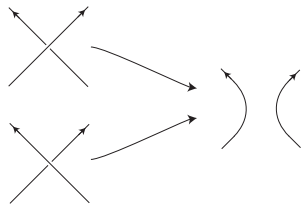}
    \caption{Nullification moves}
  \label{fig:nullification}
  \end{center}
\end{figure}

The coherent band surgeries are essentially equivalent to the nullification moves 
as shown in Figure \ref{fig:nullification}.
For an oriented link, the {\em nullification number} of the link is the minimum number of nullification moves 
needed to obtain an unlink from the link. 
This concept of nullification number has been investigated by both mathematicians and biologists, see e.g. \cite{C,S,E-M-S} and   
\cite{D-E-M}, where it is termed the \textit{general nullification number}.   (We discuss the biological applications of this in Section \ref{s:Bio}.) 
In particular, \cite{E-M-S} determined the nullification number for all prime knots up to $9$ crossings. 

For two oriented links $L$ and $L'$, we define the {\em coherent band-Gordian distance} between $L$ and $L'$, 
and denote by $d(L,L')$, as the minimum number of  coherent band surgeries needed to obtain $L$ from $L'$. 
Note that (coherent) band surgeries are reversible:  the band for a band surgery 
$L\to L'$ is also a band for a band surgery $L'\to L$, 
and so nullification moves are reversible as well.   

\subsection{Current Work.}

Our present interest lies in coherent band surgeries. We consider both the minimal (coherent band surgery or nullification) `distance' between knots and links, and the \textit{pathways} between knots and multi-component links.  In particular, we examine when there are $n$-component link intermediates (for $n \geq 1$) between two given knots or links.  

Our main results are to categorise coherent band pathways between knots and multi-component links in Section \ref{s:CBP}.  
Additionally, we characterise the minimal coherent band pathway (with intermediates) between any two knots or 2-component links with small crossing number.  We demonstrate these band surgeries for knots and links with small crossing numbers in Section \ref{s:Tables}.  In Tables \ref{Table:KK}, \ref{Table:KL} and \ref{Table:LL}, we also demonstrate an intermediate in this minimal pathway, although note in general there could be
more than one pathway of minimal length (and hence several possible intermediates). 
For these minimal pathways, we also demonstrate where the band occurs in Figures  \ref{fig:bandexamples-5} and \ref{fig:bandexamples6-7}.

In Section \ref{s:Bio}, we apply these results  to give lower bounds 
on the minimum number of recombinant events separating DNA configurations, restrict the recombination pathways and determine
chirality and/or orientation of the resulting recombinant DNA molecules.

\subsection{Acknowledgements}
The authors would like to thank Claus Ernst, Akio Kawauchi, Taizo Kanenobu and Robert Scharein for their helpful discussions.  
This work is partially supported by The Leverhulme Trust Research Grant RP2013-K-017 and EPSRC grants G039585/1 and H031367/1 to Dorothy Buck, and EPSRC grant H031367/1 and the Japan Society for the Promotion of Science KAKENHI 26800081 to Kai Ishihara.  
\section{Coherent Band Pathways}
\label{s:CBP}
A sequence $L_0,L_1,\ldots ,L_n$ of oriented links is called a 
{\em coherent band pathway with length} $n$, denoted by $L_0\leftrightarrow L_1\leftrightarrow \ldots \leftrightarrow L_n$, if $d(L_{i-1},L_i)=1$ for any integer 
$i\in \{\,1,2,\ldots, n\}$.
The number of components of a link $L$ is denoted by $\mu(L)$. 
A coherent band surgery  always changes 
$\mu(L)$ by one, thus:  
\begin{rem}\label{rem:comp}
$d(L,L')\equiv \mu(L)-\mu(L') \pmod{2}$.
\end{rem}
We consider the number of link components on a coherent band pathway.
Let $L_0\leftrightarrow L_1\leftrightarrow L_2\leftrightarrow L_3$ 
be a coherent band pathway with length $3$, then we have the following:
\begin{lem}\label{lem:comp}
$(1)$ If the links $L_0,L_1$ and $L_2$ have $m,m-1$ and $m$ components respectively,
then there exists an $(m+1)$-component link $L'_1$ such that 
$L_0\leftrightarrow L'_1\leftrightarrow L_2$ 
is a coherent band pathway with length $2$ 
(see  Figure \ref{fig:compchange1}, left).\\
$(2)$ If the links $L_0,L_1,L_2$ and $L_3$ have $m-1,m,m+1$ and $m$ components respectively, 
then there exist an $m$-component link $L'_1$ and $(m-1)$-component link $L'_2$  
such that 
$L_0\leftrightarrow L'_1\leftrightarrow L'_2\leftrightarrow L_3$ is a coherent band pathway with length $3$ 
(see Figure \ref{fig:compchange1}, right).
\end{lem}

\begin{figure}[htbp]
\begin{center}
\mbox{}\hfill
\begin{picture}(450,50)
 \put(0,0){$(m-1)$-component} 
 \put(20,20){$m$-component} 
 \put(0,40){$(m+1)$-component} 
 \put(100,20){$L_0$} 
 \put(120,0){$L_1$}
 \put(140,20){$L_2$}
 \put(155,20){$\Rightarrow$}
 \put(170,20){$L_0$} 
 \put(190,40){$L'_1$}
 \put(210,20){$L_2$}
 \put(110,13){\rotatebox{-45}{$\leftrightarrow$}}
 \put(130,7){\rotatebox{45}{$\leftrightarrow$}}
 \put(180,27){\rotatebox{45}{$\leftrightarrow$}}
 \put(200,33){\rotatebox{-45}{$\leftrightarrow$}}
 \put(250,0){$L_0$} 
 \put(270,20){$L_1$}
 \put(290,40){$L_2$}
 \put(310,20){$L_3$} 
 \put(325,20){$\Rightarrow$}
 \put(340,0){$L_0$} 
 \put(360,20){$L'_1$}
 \put(380,0){$L'_2$}
 \put(400,20){$L_3$} 
 \put(260,7){\rotatebox{45}{$\leftrightarrow$}}
 \put(280,27){\rotatebox{45}{$\leftrightarrow$}}
 \put(300,33){\rotatebox{-45}{$\leftrightarrow$}}
 \put(350,7){\rotatebox{45}{$\leftrightarrow$}}
 \put(370,13){\rotatebox{-45}{$\leftrightarrow$}}
 \put(390,7){\rotatebox{45}{$\leftrightarrow$}}
\end{picture}
\hfill\mbox{}
\caption{}
\label{fig:compchange1}
\end{center}
\end{figure}

\begin{figure}[htbp]
  \begin{center}
    \includegraphics[scale=0.6]{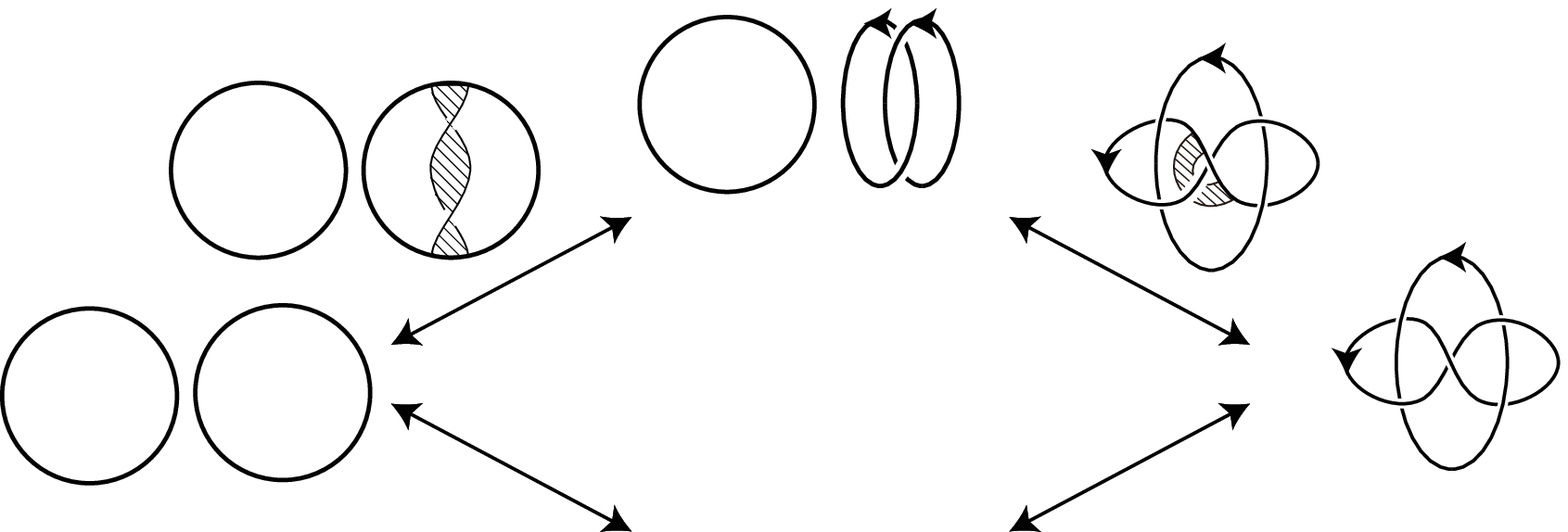}

\begin{picture}(0,0)(200,0)
\put(180,0){No knots}
\put(185,65){$0_1\sqcup2_1^2$}
\put(85,10){$0_1^2$}
\put(320,10){$5_1^2$}
\end{picture}
    \caption{}
  \label{fig:unlink-whitehead}
  \end{center}
\end{figure}

\begin{rem}\label{rem:comp}
The converse of Lemma  \ref{lem:comp},(1) is not true in general. 
In fact, there are examples of pairs of two components links which are coherent band-Gordian distance $2$, but there are no knots as an intermediate, see Figure \ref{fig:unlink-whitehead} and also Table \ref{Table:LL}.
\end{rem}

\begin{proof}
By the definition of a coherent band pathway, there exist three bands $b_1,b_2$ and $b_3$ 
such that $b_i$ relates $L_{i-1}$ and $L_{i}$, 
i.e. there exists two disjoint unions $\alpha_i\subset L_{i-1}$ and $\beta_i \subset L_{i}$ of two arcs such that
$\partial b_i=\alpha_i\cup \beta_i$, 
$\alpha_i\cap \beta_i=\partial \alpha_i=\partial \beta_i$, 
$L_{i-1}\cap b_i=\alpha_i$, $L_{i}\cap b_i=\beta_i$,
and $L_i=(L_{i-1}-\alpha_{i})\cup \beta_{i}$
($L_{i-1}=(L_{i}-\beta_{i})\cup \alpha_{i}$) for each $i\in\{1,2,3\}$. 
We may assume three bands $b_1,b_2,b_3$ are mutually disjoint 
by moving them along the links if necessary. 
Each $b_i$ is a band of a coherent band surgery, so $\mu(L_i)-\mu(L_{i-1})=1$
 or $-1$ according to two arcs of $\alpha_i$ are situated 
 in the same component of $L_{i-1}$ or situated in different components.

In the case (1) of Lemma \ref{lem:comp}, two arcs of $\alpha_1$ 
(resp. $\alpha_2$) are situated in different components of $L_0$ 
(resp. situated in the same component of $L_1$).
Let $L'_1$ be a link obtained from $L_0$ by a coherent band surgery 
along $b_2$, i.e. $L'_1=(L_0-\alpha_2)\cup \beta_2=(L_2-\beta_1)\cup\alpha_1$. Note that $L_0\leftrightarrow L'_1\leftrightarrow L_2$ is a coherent band pathway. We will show that $b_2$ can be moved along $L_1$ so that $\mu(L'_1)-\mu(L_0)=1$ i.e. $\mu(L'_1)=m+1$.  
If two arcs of $\alpha _{2}$ are situated in the same component of $L_0$, then it is not necessary to move $b_2$. 
Therefore we suppose two arcs of $\alpha _{2}$ are situated in different components of $L_0$. Then the two components agree with ones in which two arcs of $\alpha_{1}$ are situated. Because two arcs of $\alpha_2$ are situated in the same component of $L_1=(L_0-\alpha_1)\cup \beta_1$, while they are situated in different components of $L_0$. It implies that each arc of $\beta_1$ connects the two components. Hence we can move one arc of $\alpha_2$ together with $b_2$ 
along $L_0$ through one arc of $\beta_1$. Then $\alpha _{2}$ are situated in the same component of $L_0$ after that, and so $\mu(L'_1)=m+1$. 
   
In the case (2) of Lemma  \ref{lem:comp}, two arcs of $\alpha_1$ 
(resp. $\alpha_2$) are situated in the same component of $L_0$ (resp. $L_1$) and $\alpha_3$ are situated in different components of $L_2$.
If $\alpha_3$ are situated in different components of $L_1$, 
by putting $L'_1:=L_1$ and $L'_2:=(L_1-\alpha_3)\cup \beta_3=(L_3-\beta_2)\cup\alpha_2$, then we have a coherent band pathway $L_0\leftrightarrow L'_1\leftrightarrow L'_2\leftrightarrow L_3$ with $\mu(L'_1)=m$ and $\mu(L'_2)=m-1$. 
Therefore we suppose $\alpha_3$ are situated in the same component of $L_1$. The component agrees with the one which contains both arcs of $\alpha_2$.  Because two arcs of $\alpha_3$ are situated in different components of $L_2=(L_1-\alpha_2)\cup \beta_2$, while they are situated in the same component of $L_1$. 
Let $K_1$ be such a component of $L_1$ containing both $\alpha_2$ and $\alpha_3$, and for $i\in\{2,3\}$, let $\alpha_{i,1},\alpha_{i,2}$ be two arcs of $\alpha_i$.  One of arcs $\alpha_{2,1},\alpha_{2,2}$ and one of $\alpha_{3,1},\alpha_{3,2}$ appear alternately as one goes along $K_1$. say in order by $\alpha_{2,1}, \alpha_{3,1}, \alpha_{2,2}$ and $\alpha_{3,2}$, otherwise $\alpha_{3,1}$ and  $\alpha_{3,2}$ are situated in the same component of $L_2$. 
Since at most one arc of $\beta_1$ lie in $K_1$, 
one component of $L_0=(L_1-\beta_1)\cup \alpha_1$, say $K_0$, contains 
both $\alpha_2$ and $\alpha_3$. Moreover four arcs also appear in order by $\alpha_{2,1}, \alpha_{3,1}, \alpha_{2,2}$ and $\alpha_{3,2}$ as one goes along $K_0$. 
Then by putting $L'_1:=(L_0-\alpha_2)\cup \beta_2$ and 
$L'_2:=(L'_1-\alpha_3)\cup\beta_3=(L_3-\beta_1)\cup \alpha_1$,
we have a coherent band pathway
$L_0\leftrightarrow L'_1\leftrightarrow L'_2\leftrightarrow L_3$ with $\mu(L'_1)=m$ and $\mu(L'_2)=m-1$. 
\end{proof}

By using Lemma \ref{lem:comp}, we have the following.
\begin{prop}\label{prop:comp}
$(1)$ Let $L,L'$ be $m$-component links.
If $d(L,L')\ge 4$, then there exists an $m$-component link $L''$
such that $d(L,L')=d(L,L'')+d(L'',L')$ and $2\le d(L,L'')\le d(L,L')-2$.\\
$(2)$ Let $L$ $L'$ be an $m$-component link and $m'$-component link respectively. 
Then there exists a $m$-component link $L''$ such that $d(L,L')=d(L,L'')+d(L'',L')$ and $d(L'',L')=|m-m'|$.
\end{prop}

By the same argument of Baader \cite{B}, we have the following.
\begin{thm}\label{thm:inftylinks}
Let $L$ and $L'$ be oriented links with the same number of components.  
If $d(L,L')=2$, then there exists infinite family $\{M_{i}\}$ of mutually distinct links such that $d(L,M_{i})=d(M_{i},L')=1$; i.e. $L\leftrightarrow M_{i}\leftrightarrow L'$ is a coherent band pathway with length $2$.
\end{thm}
\begin{proof}
Since $d(L,L')=2$, there exists an oriented link $M$ with $d(M,L)=d(M,L')=1$. We may assume there exist two disjoint bands $b$ and $b'$ such that 
$M\cap b=\alpha\subset\partial b$, 
$M\cap b'=\alpha\subset\partial b'$, and $L=(M-\alpha)\cup\beta$, $L'=(M-\alpha')\cup\beta'$, where each of $\alpha,\beta,\alpha',\beta'$ is a disjoint union of two arcs such that $\partial b=\alpha\cup\beta$, $\alpha\cap\beta=\partial\alpha=\partial\beta$. $\partial b'=\alpha'\cup\beta'$, $\alpha'\cap\beta'=\partial\alpha'=\partial\beta'$. 
There exist eight types for the pair $(M, b\cup b')$ according to the situation of four arcs $\alpha\cup\alpha'$ in $M$, see Figure \ref{fig:length2}
\begin{figure}[htbp]
  \begin{center}
    \includegraphics[scale=0.9]{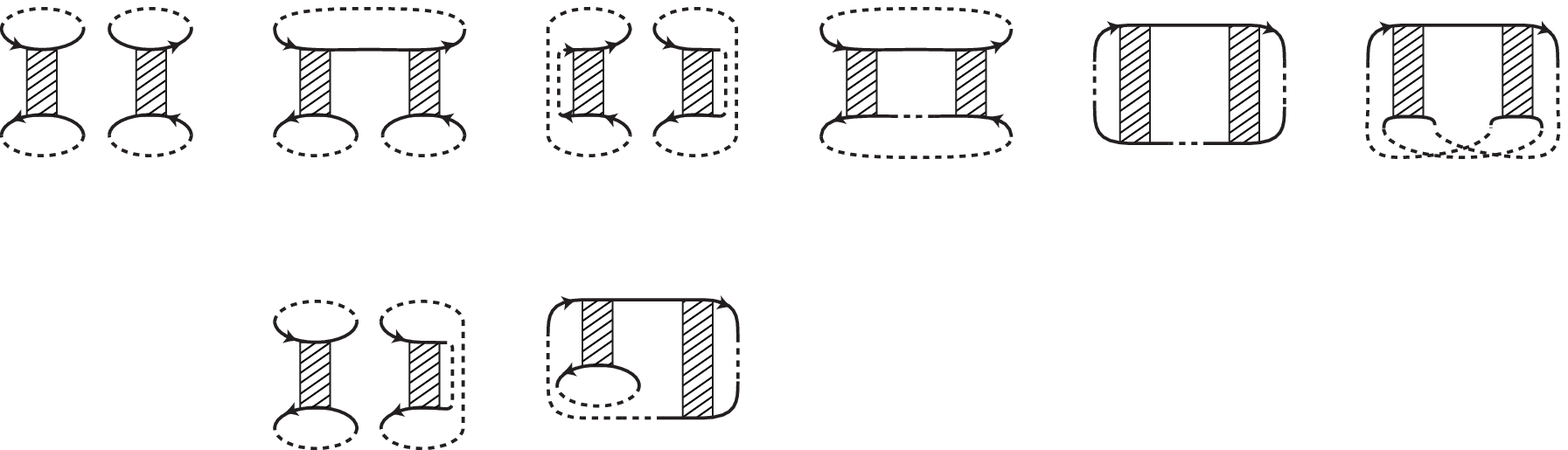}

\begin{picture}(500,0)
\put(5,145){$3\leftrightarrow 4\leftrightarrow3$}
\put(85,145){$2\leftrightarrow 3\leftrightarrow2$}
\put(165,145){$3\leftrightarrow 2\leftrightarrow3$}
\put(250,145){$1\leftrightarrow 2\leftrightarrow1$}
\put(327,145){$2\leftrightarrow 1\leftrightarrow2$}
\put(410,145){$2\leftrightarrow 1\leftrightarrow2$}
\put(85,60){$2\leftrightarrow 3\leftrightarrow4$}
\put(165,60){$1\leftrightarrow 2\leftrightarrow3$}
\put(20,90){(1)}
\put(100,90){(2-1)}
\put(180,90){(3-1)}
\put(260,90){(3-2)}
\put(340,90){(4-1)}
\put(423,90){(4-2)}
\put(100,0){(2-2)}
\put(180,0){(3-3)}

\end{picture}

    \caption{}
  \label{fig:length2}
  \end{center}
\end{figure}

Let $M^*$ be an oriented link obtained from $L$ (resp. $L'$) by 
a coherent band surgery along a band $b'$ (resp. $b$), i.e $M^*=(L-\alpha')\cup\beta'=(L'-\alpha)\cup\beta$.
Note that $M^*\cup b\cup b'=M\cup b\cup b'$ and $d(M^*,L)=d(M^*,L')=1$, so we call the pair $(M^*, b\cup b')$ the {\em dual} of $(M,  b\cup b')$. 
It is enough to prove the statement of Theorem \ref{thm:inftylinks} for the case where each component of $M$ intersects with $b$ or $b'$. Since  $M\cap (b\cup b')=\alpha\cup\alpha'$ is a union of mutually disjoint four arcs, we may assume the number of components $\mu(M)$ is at most $4$. There are four cases: (1) $M$ has $4$-components; (2) $M$ has $3$-components; 
(3) $M$ has $2$-components; (4) $M$ is a knot.

In the first case (Figure \ref{fig:length2}, (1)), the dual of $(M^*,b\cup b')$ is a pair of $2$-component link and two bands (Figure \ref{fig:length2}, (3-1)). Hence this reduces to the third case. 

In the second case, 
one component contains one arc of $\alpha$ and one arc of $\alpha'$, and other two components each contain only one arc of $\alpha\cup\alpha'$ (Figure \ref{fig:length2}, (2-1)), otherwise $L$ and $L'$ have different numbers of components (Figure \ref{fig:length2}, (2-2)). In this case, the dual of $(M^*,b\cup b')$ is a pair of a knot and two bands (Figure \ref{fig:length2}, (4-1)). Hence this reduces to the fourth case. 

In the third case, 
either each component contains one of $\alpha$ and $\alpha'$ 
(Figure \ref{fig:length2}, (3-1)) or each component intersects 
both $\alpha$ and $\alpha'$ (Figure \ref{fig:length2}, (3-2)), 
otherwise $L$ and $L'$ have different numbers of components 
(Figure \ref{fig:length2}, (3-3)). 
Let $\beta_n,\beta'_n$ be unions of two arcs obtained 
from $\beta,\beta'$ by $n$ time full twisting each other 
together with bands $b,b'$, and put $M_n:=(M-\alpha)\cup\beta_n$ 
as shown in Figure \ref{fig:inftylinks}, (3-1) and (3-2).  
Note that $L\leftrightarrow M_n\leftrightarrow L'$ is a coherent band pathway.
\begin{figure}[htbp]
  \begin{center}
    \includegraphics[scale=0.9]{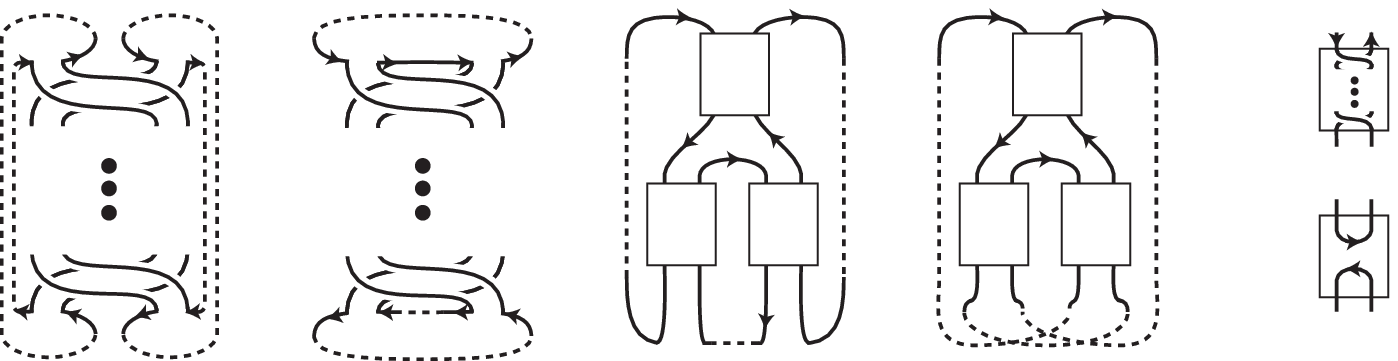}

\begin{picture}(400,0)
\put(35,0){(3-1)}
\put(118,0){(3-2)}
\put(198,0){(4-1)}
\put(280,0){(4-2)}
\put(207,85){$x$}
\put(193,45){$y$}
\put(220,45){$z$}
\put(288,85){$x$}
\put(275,45){$y$}
\put(300,45){$z$}
\put(380,80){$\biggr\}2n$}
\put(335,72){\framebox(15,22){$n$}}
\put(352,80){$=$}
\put(335,30){\framebox(15,22){$\infty$}}
\put(352,38){$=$}
\end{picture}

    \caption{}
  \label{fig:inftylinks}
  \end{center}
\end{figure}
In the case (3-1), $M_n$ has four components $K_{n,1},K_{n,2},K_{n,3},$ 
and $K_{n,4}$, 
where $K_{n,1}\cup K_{n,2}$ contains $\beta_n$ and $K_{n,3}\cup K_{n,4}$ contains $\beta'_n$.
Then we have the following. 
\begin{align*}(lk(K_{n,1},K_{n,3}),lk(K_{n,1},K_{n,4}),lk(K_{n,2},K_{n,3}),lk(K_{n,2},K_{n,4}))\\
=(lk(K_{0,1},K_{0,3})\pm n,lk(K_{0,1},K_{0,4})\mp n,lk(K_{0,2},K_{0,3})\mp n,lk(K_{0,2},K_{0,4})\pm n)
\end{align*}
Hence the family $\{M_n\}_{n\in\mathbb{Z}}$ contains infinitely many different links. 
In the case (3-2), $M_n$ has two components $K_{n,1}$ and $K_{n,2}$.
Then we have $lk(K_{n,1},K_{n,2})=lk(K_{0,1},K_{0,2})\pm 2n$. Hence the family $\{M_n\}_{n\in\mathbb{Z}}$ contains infinitely many different links. 

In the fourth case (Figure \ref{fig:length2}, (4-1) and (4-2)), we will use the same argument as \cite{B}. 
For $x,y,z\in \mathbb{Z}\cup\{\infty\}$, $M_{x,y,z}$ denotes the oriented link as shown in Figure \ref{fig:inftylinks}, (4-1) and (4-2), 
and put $M_n:=M_{-n,n,n}$. 
Then by Lemma 1.3 in \cite{B}, 
$$\lim_{n\to\infty}\frac{1}{n}\nabla_{M_{n}}(z)=-z^2\nabla_{M_{\infty,0,\infty}}(z),$$ where $\nabla$ is the Conway polynomial.  
Since $M_{\infty,0,\infty}$ is a knot, which is ambient isotopic to $M$, $\nabla_{M_{\infty,0,\infty}}(z)$ is not zero. 
Hence the family $\{M_n\}_{n\in\mathbb{Z}}$ contains infinitely many different links. 
\end{proof}

\section{Coherent Band pathways between knots and 2-component links}

\subsection{Previous lower bounds for the length of the Coherent Band Pathway}
We begin by recalling known lower bounds for coherent band pathway.
Let $L$ and $L'$ be oriented links. The signature of $L$ is denoted by $\sigma(L)$. 
Murasugi \cite{M} showed the following:
\begin{thm}[\cite{M}]\label{thm:signature}
\[|\sigma(L)-\sigma(L')|\le d(L,L').\]
\end{thm}

We denote the Jones polynomial of $L$ by $V(L; t)$ and put $\omega=e^{i\pi/3}$.
Kanenobu \cite{Kn1} showed the following.
\begin{thm}[\cite{Kn1}]\label{thm:Jones}
Suppose $\mu(L)\equiv \mu(L')\pmod{2}$. 
\begin{center}
If $V(L;\omega)/V(L';\omega)\notin\{\pm (i\sqrt{3})^{\pm k}, \sqrt{3}^{\pm n}\ |\ k=0,1,\ldots,n-1\}$, 
then $n+2\le d(L,L')$.
\end{center}
Suppose $\mu(L)\not\equiv \mu(L')\pmod{2}$. 
\begin{center}
If $V(L;\omega)/V(L';\omega)\notin\{\pm i (i\sqrt{3})^{\pm k}, -\sqrt{3}^{\pm n}\ |\ k=0,1,\ldots,n-1\}$, 
then $n+2\le d(L,L')$.
\end{center}
\end{thm}
We denote the Q polynomial of a link $L$ by $Q(L; z)$, and put $\rho(L)=Q(L; (\sqrt{5}-1)/2)$.
Kanenobu \cite{Kn1} showed the following.
\begin{thm}[\cite{Kn1}]\label{thm:Q}
If $\rho(L)/\rho(L')\notin\{\pm \sqrt{5}^{\pm k}, \sqrt{5}^{\pm n}\ |\ k=0,1,\ldots,n-1\}$, 
then $n+1\le d(L,L')$.
\end{thm}

Since a coherent band move changes the number of link components by one, 
any oriented links $L$ and $L'$ satisfy $d(L,L')\ge |\mu(L)-\mu(L')|$.
From a property of the Arf invariant, we obtain the following (see \cite{R} and \cite[Corollary 8.3.2]{Kw0}).

\begin{thm}\label{thm:Arf}
Let $L$ be a oriented proper link. 
If $d(L,L')=\mu(L)-\mu(L')$, then $L'$ is also proper and ${\rm Arf}(L)={\rm Arf}(L')$. 
In particular, if both $L$ and $L'$ are proper and $d(L,L')=1$, then
${\rm Arf}(L)={\rm Arf}(L')$.
\end{thm}

We denote the Alexandar polynomial of a knot $K$ by $\Delta(K)$. 
Kawauchi \cite{Kw} showed the following.
\begin{thm} [\cite{Kw}]\label{thm:Alex} 
Let $T'_{2,2k}$ be the anti-parallel $(2,2k)$-torus link.
If $d(K,T'_{2,2k})=1$, then there exists a polynomial $f$ such that
$$\Delta_{K}(t)\equiv \pm t^{r}f(t)f(t^{-1})\pmod{k}.$$
\end{thm}

Fox and Milnor \cite{F-M} showed a condition of the Alexandar polynomial 
$\Delta_K(t)$ for $d(K,0^2_1)=1$, 
and Kawauchi showed a condition of the Alexandar polynomial 
$\Delta_K(t)$ for $d(K,T_{2,2k})=1$ (see \cite{I-S}).
Here $T_{2,2k}$ is the parallel $(2,2k)$-torus link. 
\cite{D-I-M-S} showed a necessary and sufficient condition for $d(S(4mn-1,2m),T'_{2,2k})=1$, where $S(4mn-1,2m)$ is a $1$ genus $2$ bridge knot.

\subsection{Our Results on Minimal Coherent Band Pathways.}
\label{s:Tables}
In this section, we characterise the minimal coherent band pathways between $L$ and $L'$, where $\{L,L'\} = \{$knot, 2-component link$\}$.  
We list as illustration one intermediate in the pathway, although note in general there could be more than one pathway of minimal length.  

For unoriented knots and links, we use the Rolfsen notation (Appendix C of \cite{Rolf}.)    Given a knot or link $L$ then we write $L!$ to denote its mirror image.  For oriented 2-component links, we use the same notation as Kanenobu  \cite{Kn2}.   Namely,  we choose $L$ to be the 2-component link with negative linking number, and $L!$ to be the link with positive linking number.   Given an oriented 2-component link then we write $L'$ to denote the same link with the orientation reversed on one component.  

\subsubsection{Coherent Band Pathway between two knots}

In this section, we discuss the coherent band pathways between two knots $K$ and $K'$, and show Table \ref{Table:KK} for knots of 7 crossings or less. 
From Remark \ref{rem:comp}, $d(K,K')$ is even.
By the definition of coherent band pathway, 
$d(K,K')=0$ if and only if $K=K'$. We indicate this 
with ``$0$'' in Table \ref{Table:KK}. 
$d(K,K')=2$ if and only if $K\neq K'$ and there exists a ($2$-component)
 link $L$ such that $d(K,L)=d(L,K')=1$. We indicate this 
with ``$2(L)$'' in Table \ref{Table:KK}.  
In the case where $d(K,K')=d\ge 4$, we need to use the lower bounds,  
Theorems \ref{thm:signature}, \ref{thm:Jones}, \ref{thm:Q}, 
and the fact that there exists a knot $K''$ such that $d(K,K')=d(K,K'')+d(K'',K')$ 
and $2\le d(K,K'')\le d(K,K')-2$ from Proposition \ref{prop:comp} (1).
 We indicate this with ``$d($I$K'')$'', ``$d($II$K'')$'',``$d($III$K'')$'' 
respectively in Table \ref{Table:KK}.

The second author presented an earlier version of this table (for knots of 7 crossings or less) 
at the Mathematics of Knots in Tokyo Workshop in December 2011.  In the subsequent conference proceedings article \cite{us-Japan}, the corresponding 
table considers knots of 6 crossings or less because of space constraints.   We have since learned that Kanenobu and Moriuchi \cite{Kn-Mo} 
have also independently extended our table in \cite{us-Japan} to include knots with 7 crossings.   They also corrected an entry (to 4, between the knots $3_1$ and $6_1$) in our earlier version.  

\begin{table}
\caption{Coherent Band Pathways between two knots}\label{Table:KK}
\begin{center}
\scalebox{0.9}{
\input{DistanceKK1.tex}

}
\end{center}
\end{table}

\begin{table}
{\scshape Table 1. }(Continued)
\begin{center}
\scalebox{0.9}{
\input{DistanceKK2.tex}

}
\end{center}
\begin{flushleft}
{\small
$\dagger$: Corrected by Kanenobu and Moriuchi \cite{Kn-Mo}.

$\dagger\dagger$: Improved by Kanenobu and Moriuchi \cite{Kn-Mo}.
}
\end{flushleft}
\end{table}
\subsubsection{Coherent Band Pathway between knots and 2-component links}
\begin{figure}[htbp]
  \begin{center}
    \includegraphics[scale=0.6]{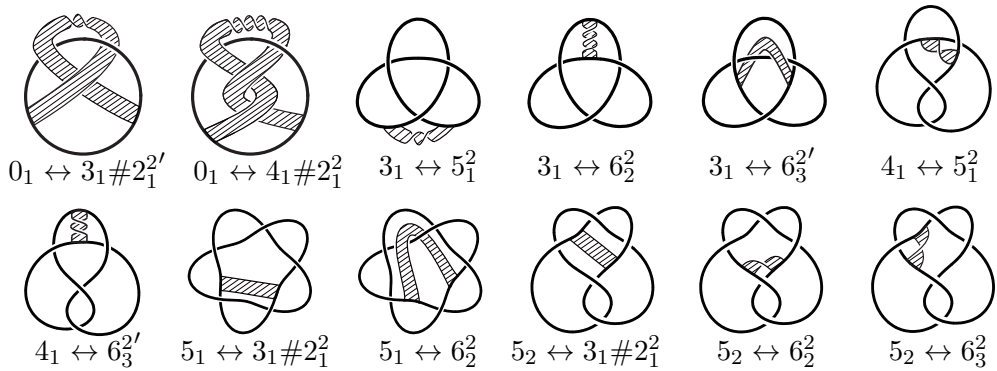}

\begin{picture}(0,0)(200,0)
\put(10,73){$0_1\leftrightarrow 3_1\#{2_1^2}' $}
\put(80,73){$0_1\leftrightarrow 4_1\#2_1^2$}
\put(150,75){$3_1\leftrightarrow 5_1^2$}
\put(210,75){$3_1\leftrightarrow 6_2^2$}
\put(275,75){$3_1\leftrightarrow {6_3^2}'$}
\put(340,75){$4_1\leftrightarrow 5_1^2$}

\put(20,5){$4_1\leftrightarrow {6_3^2}' $}
\put(75,5){$5_1\leftrightarrow 3_1\#2_1^2$}
\put(150,5){$5_1\leftrightarrow 6_2^2$}
\put(200,5){$5_2\leftrightarrow 3_1\#2_1^2$}
\put(278,5){$5_2\leftrightarrow 6_2^2$}
\put(343,5){$5_2\leftrightarrow 6_3^2$}

\end{picture}
    \caption{Bands attaching to knots with up to $5$ crossings}
  \label{fig:bandexamples-5}
  \end{center}
\end{figure}

\begin{figure}[htbp]
  \begin{center}
    \includegraphics[scale=0.6]{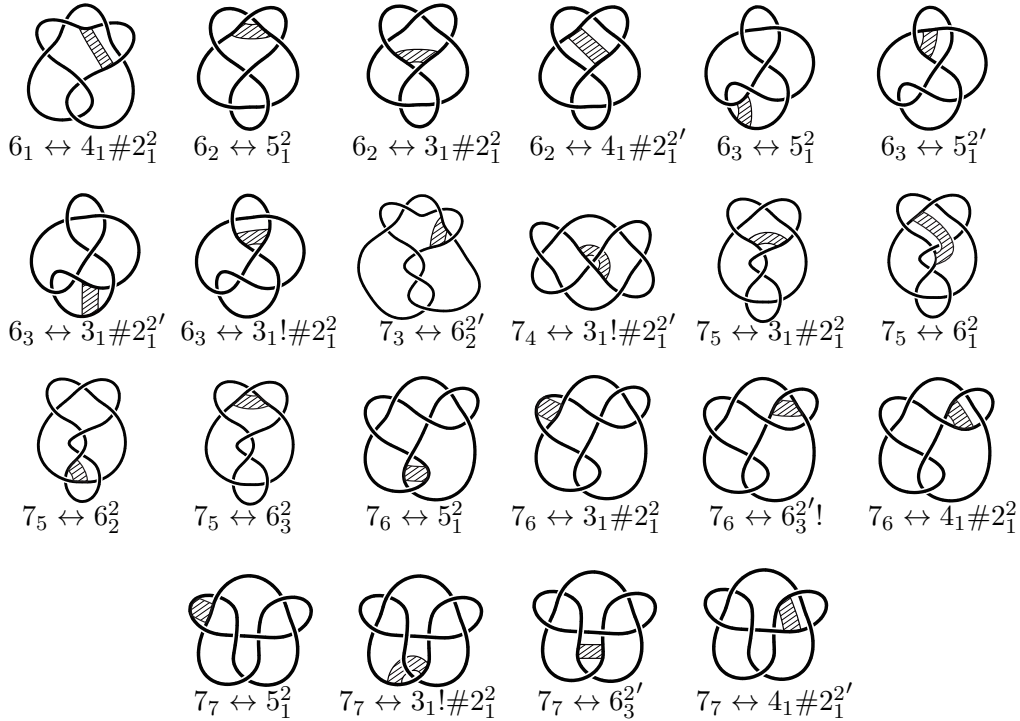}

\begin{picture}(0,0)(200,0)
\put(10,215){$6_1\leftrightarrow 4_1\#{2_1^2} $}
\put(80,215){$6_2\leftrightarrow 5_1^2$}
\put(140,215){$6_2\leftrightarrow 3_1\#2_1^2$}
\put(207,215){$6_2\leftrightarrow 4_1\#{2_1^2}' $}
\put(278,215){$6_3\leftrightarrow 5_1^2$}
\put(340,215){$6_3\leftrightarrow {5_1^2}'$}

\put(10,145){$6_3\leftrightarrow 3_1\#{2_1^2}' $}
\put(75,145){$6_3\leftrightarrow 3_1!\#{2_1^2}$}
\put(150,145){$7_3\leftrightarrow {6_2^2}'$}
\put(200,145){$7_4\leftrightarrow 3_1!\#{2_1^2}' $}
\put(270,145){$7_5\leftrightarrow 3_1\#{2_1^2}$}
\put(340,145){$7_5\leftrightarrow {6_1^2}$}

\put(15,75){$7_5\leftrightarrow {6_2^2} $}
\put(80,75){$7_5\leftrightarrow {6_3^2}$}
\put(145,75){$7_6\leftrightarrow {5_1^2}$}
\put(200,75){$7_6\leftrightarrow 3_1\#{2_1^2} $}
\put(275,75){$7_6\leftrightarrow {6_3^2}'!$}
\put(335,75){$7_6\leftrightarrow 4_1\#{2_1^2}$}

\put(80,5){$7_7\leftrightarrow {5_1^2}$}
\put(135,5){$7_7\leftrightarrow 3_1!\#{2_1^2}$}
\put(210,5){$7_7\leftrightarrow {6_3^2}' $}
\put(270,5){$7_7\leftrightarrow 4_1\#{2_1^2}'$}
\end{picture}
    \caption{Bands attaching to knots with $6$ or $7$ crossings}
  \label{fig:bandexamples6-7}
  \end{center}
\end{figure}

In this section, we discuss the coherent band pathway between a knots $K$ and a $2$-component link $L$, and show Table \ref{Table:KL}. 
From Remark \ref{rem:comp}, $d(K,L)$ is odd.
By the definition of coherent band pathway, 
$d(K,K')=1$ if and only if there exists a coherent band surgery between $K$ and $L$. 
Examples of band surgeries are shown in Figure \ref{fig:bandexamples-5}, \ref{fig:bandexamples6-7} and \cite{I-S}. 
We indicate this with ``$1$'' in Table \ref{Table:KL}. 
In the case where $d(K,L)=d\ge 3$, we need to use the lower bounds,  
Theorems \ref{thm:signature}, \ref{thm:Jones}, \ref{thm:Q}, \ref{thm:Arf}, \ref{thm:Alex}, and the fact that there exists a knot $K'$ such that $d(K,L)=d(K,K')+d(K',L)$
and $d(K',L)=1$ from Proposition \ref{prop:comp} (2).
 We indicate this with ``$d($I$K')$'', ``$d($II$K')$'', 
 ``$d($III$K')$'',  ``$d($IV$K')$'', ``$d($V$K')$''
respectively in Table \ref{Table:KL}.

\begin{table}
\caption{Coherent band pathways between knots and 2-component links}\label{Table:KL}
\begin{center}
\scalebox{0.85}{
\input{DistanceKL1.tex}

}
\end{center}
\end{table}

\begin{table}
{\scshape Table 2. }(Continued)
\begin{center}
\scalebox{0.9}{
\input{DistanceKL2.tex}

}
\end{center}
\begin{flushleft}
{\small 
$\dagger$: Corrected by Kanenobu and Moriuchi \cite{Kn-Mo}.

$\dagger\dagger$: Improved by Kanenobu and Moriuchi \cite{Kn-Mo}.

$\dagger\dagger\dagger$: Pointed out by Robert Scharein and also independently improved by Kanenobu and Moriuchi \cite{Kn-Mo}.
}
\end{flushleft}
\end{table}
\subsubsection{Coherent Band Pathways between two 2-component links}
In this section, we discuss the coherent band pathway between two 
$2$-component links $L$ and $L'$, and show Table \ref{Table:LL}. 
From Remark \ref{rem:comp}, $d(L,L')$ is even.
By the definition of coherent band pathway, 
$d(L,L')=0$ if and only if $L=L'$. We indicate this 
with ``$0$'' in Table \ref{Table:LL}. 
If $d(L,L')=2$, there are two cases: (1) there exists a knot $K$ 
such that $d(L,K)=d(K,L')=1$ (we indicate this with ``$2(K)$''), 
(2) there exists a $3$-component link $L''$ such that 
$d(L,L'')=d(L'',L')=1$ (we indicate this with ``$2(L'')$'').   
In the case (1), by Lemma \ref{lem:comp} (1), the case (2) happens as well.
In the case (2), however, the case (1) is not necessary to happen. 
In fact,  by Theorem \ref{thm:Arf}, the case (2) can happen but (1) if $L$ and $L'$ are proper and Arf$(L)\neq$Arf$(L')$. 
We indicate this with ``*'' in Table \ref{Table:LL}. 

\begin{thm}\label{thm:2comp}
Let $L$ and $L'$ be $2$-component links. 
Suppose there is a $3$-component link $L''$ such that $d(L,L'')=d(L'',L')=1$ (so $d(L,L')\le d(L,L'')+d(L'',L')=2$),
however $d(L,K)+d(K,L')> 2$ for any knot $K$. 
Then $L$ and $L'$ have the same knot as a component, 
and \[lk(L)=lk(L').\]
\end{thm}

\begin{proof}
By the same argument as the proof of Lemma \ref{lem:comp}, we may assume that 
there exist two disjoint bands $b_1$ and $b_2$ such that $b_1$ relates $L$ and $L''$, and $b_2$ relates $L''$ and $L'$. 
Let $L=K_1\cup K_2$. 
Since $L''$ is a $3$-component link, $b_1$ is attaching to one component of $L$, say $K_1$. 
In other words, $b_1$ is disjoint from $K_2$. 
Then $b_2$ is also disjoint from $K_2$, and so $L'$ has also $K_2$ as a component.
Otherwise the sequence of two coherent band surgeries along $b_2$ and $b_1$ in this order gives a coherent band pathway $L\leftrightarrow K\leftrightarrow L'$, where $K$ is a knot. 
Let $K_1\leftrightarrow (K_{11}\cup K_{12})\leftrightarrow K_1'$ be the coherent band pathway given by $b_1$ and $b_2$, i.e. $L''=K_{11}\cup K_{12}\cup K_2$ and $L'=K_1'\cup K_2$. 
Since $K_2$ does not intersect the band $b_1$ (resp., $b_2$), $lk(K_1,K_2)=lk(K_{11},K_2)+lk(K_{12},K_2)$
(resp., $lk(K_1',K_2)=lk(K_{11},K_2)+lk(K_{12},K_2)$). 
Hence $lk(K_1,K_2)=lk(K_1',K_2)$ ($lk(L)=lk(L')$).
\end{proof}

The following corollary follows directly from Theorem \ref{thm:Arf} 
and  Theorem \ref{thm:2comp}. 
\begin{cor}\label{cor:2comp}
Let $L$ and $L'$ be $2$-component links with even linking numbers. 
Suppose that Arf$(L)\neq $Arf$(L')$ and $lk(L)\neq lk(L')$.
Then $d(L,L')\ge4$.
\end{cor}

In the case where $d(L,L')=d\ge 4$, we need to use the lower bounds,  
Theorems \ref{thm:signature}, \ref{thm:Jones}, \ref{thm:Q}, Corollary \ref{cor:2comp}
and there exists a $2$-component link $L''$ such that $d(L,L')=d(L,L'')+d(L'',L')$ 
and $2\le d(L,L'')\le d(L,L')-2$ from Proposition \ref{prop:comp} (1). 
We indicate this with ``$d($I$L'')$'', ``$d($II$L'')$'', 
``$d($III$L'')$'', ``$d($IV$L'')$'' respectively in Table \ref{Table:LL}.

\begin{table}
\caption{Coherent band pathways between two 2-component links}\label{Table:LL}
\begin{center}
\scalebox{0.9}{
\input{DistanceLL1.tex}

}
\end{center}
\end{table}

\begin{table}
{\scshape Table 3. }(Continued)
\begin{center}
\scalebox{0.9}{
\input{DistanceLL2.tex}

}
\end{center}
\begin{flushleft}
{\small 
$\dagger\dagger$: Improved by Kanenobu and Moriuchi \cite{Kn-Mo}.

$*$: $d(L,L')=2$ since there exists a $3$-component link $L''$ with $d(L,L'')=d(L'',L')=1$. 
However there exists no knot $K$ with $d(L,K)=d(K,L')=1$ which is shown by Theorem  \ref{thm:Arf}. 
 
$**$: $d(L,L')=2$ since there exists a $3$-component link $L''$ with $d(L,L'')=d(L'',L')=1$. 
However we could not find a $K$ with $d(L,K)=d(K,L')=1$, and could not show no existance of such a knot. 
}
\end{flushleft}
\end{table}
%


\section{Biological Applications} 
\label{s:Bio}
%
%
 
\subsection{Site-specific Recombination.}
The current work here models the action of a family of proteins, \textit{site-specific recombinases}, acting on DNA molecules. If the axis of the famous DNA double helix is circular, site-specific recombinases can convert these circular DNA molecules into a variety of nontrivial knots and links.  

These proteins mediate site-specific recombination, the reshuffling of the genetic sequence, for example changing GATTACA into ACATTAG. The result of site-specific recombination is the deletion, insertion or inversion of a DNA segment. This corresponds to a wide variety of physiological processes, including crucial steps in viral infections \cite{Grindley6}. In addition to their inherent biochemical interest, the pharmaceutical and agricultural industries utilise site-specific recombinases as tools for precisely manipulating DNA.   
 
A site-specific recombinase first seeks two \textit{sites}: copies of a specific short string of basepairs in the DNA sequence of the original molecule(s) (the \textit{substrate}). When it finds these, it will bind to the DNA at each of these sequences and bring two in close proximity. After an intricate process of cutting and rejoining, the protein then releases the DNA molecule(s) (the \textit{product}). This process of rearranging the DNA sequence at a particular sequence of base pairs is called site-specific recombination.  Under favourable conditions, some recombinases will perform \textit{processive} recombination -- that is, the recombinase will bind and perform multiple rounds of rearranging before releasing the DNA molecules.  
 
On a single circular DNA molecule, these two sites can have the same or opposite orientation:  the sequences are palindromes (\textit{inverted repeats}) or direct copies (\textit{direct repeats}).  When site-specific recombination occurs on direct repeats, the product has one more (or less) component than the substrate:
\textit{e.g.} site-specific recombination produces a link from a knot or \textit{v.v.}.  When site-specific recombination occurs on inverted repeats, the number of components in the product is the same as the substrate.
 
 \subsection{Site-specific recombination and band surgery.}
In this paper, we model site-specific recombination by band surgery.  
This band surgery accurately reflects both the short length of the 
sites, and that the protein-protein interactions within the pre- and 
post-recombinant synaptic complex constrains the geometry of these sites.  (See \cite{JMB} for a review of this evidence.)
 
 A natural model for site-specific recombination on direct repeat sites is thus \textit{coherent band surgery}. 
 

\begin{figure}[htbp]
  \begin{center}
    \includegraphics[scale=0.8]{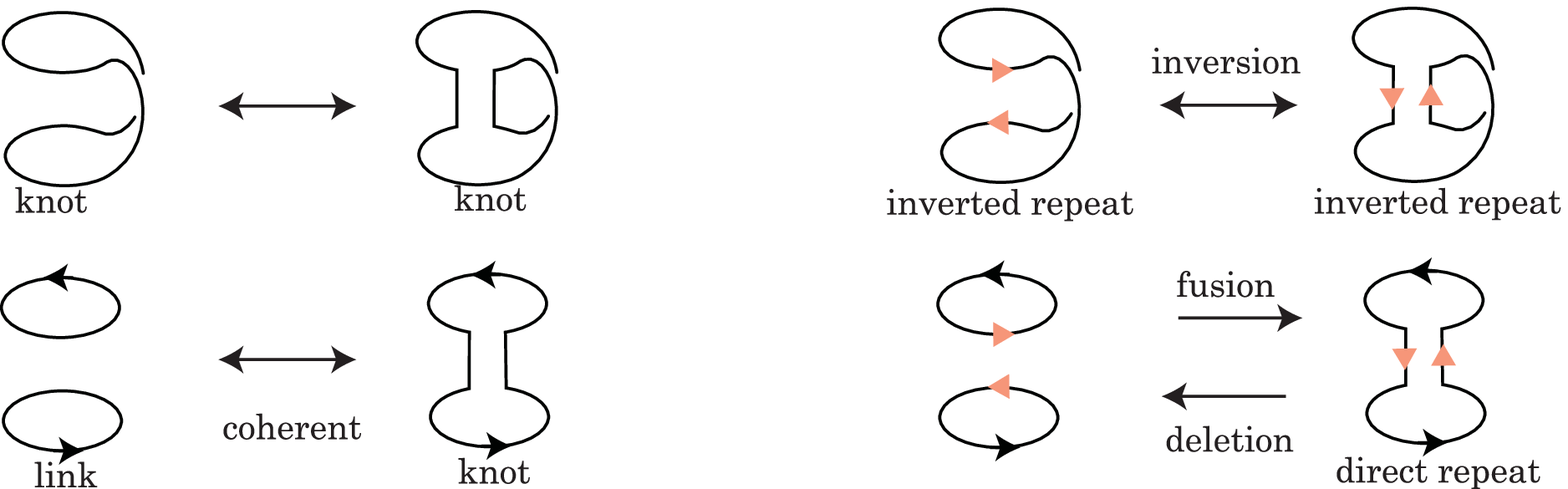}
    \caption{
    Band surgeries on the left correspond to site-specific recombinations on the right.
    In particular, a coherent band surgery on the left bottom corresponds 
    to a direct repeat recombination on the right bottom.
    }
  \label{fig:band-recombination}
  \end{center}
\end{figure}


\subsection{Examples of Biological Applications}
Although we discuss several specific biological applications below,  we emphasise that these are to illustrate the power of the the topological results in the preceding sections, and are not exhaustive.  

\subsubsection{Determining putative pathways:}
There are several ways we can articulate or restrict possible pathways of recombination, given specified products.  

For example, the site-specific recombinase Tn3 will recombine circular DNA molecules.  In \cite{Cozz1} the products of processive rounds of site-specific recombination mediated by Tn3 on an unknotted DNA substrate were shown to be $2^2_1, 4_1, 5^2_1$ and $6_2$.  An early triumph of DNA Topology was to use these knots to infer the underlying geometry of the circular DNA substrate (3 trapped supercoils), and confirm that the processive recombination pathway was of length four:  $0_1 \rightarrow 2^2_1 \rightarrow 4_1 \rightarrow 5^2_1 \rightarrow 6_2$ \cite{E-S}. 
However without the requirement of a fixed geometry there are several shorter possible processive recombination pathways, each consistent with the local strand exchange mechanism of Tn3 resolvase.  There are 3 minimal pathways between $0_1$ and $6_2$ each of length 2, with $2^2_1$, $3_1 \# 2^2_1$ and $4_1 \# 2^2_1$ as possible intermediates.  Additionally, a length 3 processive recombinantion pathway that contains most of the same products is $0_1 \rightarrow 2^2_1 \rightarrow 6_2 \rightarrow 5^2_1$, which has length 3 between the unknot and $6_2$.  

In certain circumstances the site-specific recombinase XerCD will act processively on a $6^2_1$ torus link to resolve it to the split link $0^2_1$ \cite{embo}.  Previous work of the second author and collaborators have modelled this reaction under the assumptions both that the minimal crossing number of the intermediates decreased during each stepwise reaction and that each intermediate had at most two components \cite{Shimo}.  Under these assumptions, it was shown that the most parsimonious pathway was of length 6:  \ $6^2_1 \rightarrow 5_1 \rightarrow 4^2_1 \rightarrow 3_1 \rightarrow 2^2_1 \rightarrow 0_1$, and then (no longer requiring the MCN to decrease) $0_1 \rightarrow 0^2_1$.  

Our current work demonstrates that if one relaxes either of these assumptions, there are a number of other possible pathways, all of which are consistent with the local mechanism of XerCD \cite{Colloms}.  For example, if the number of components of the intermediates could be greater than 2, then the pathway $6^2_1 \rightarrow 4^2_1 \# 2^2_1 \rightarrow 6^2_2 \rightarrow 3_1 \rightarrow 2^2_1 \rightarrow 0_1 \rightarrow 0^2_1$ is a length 6 pathway.  Similarly if one does not require that the MCN decreases at each step, then the unlinking pathway $6^2_1 \rightarrow 7_5 \rightarrow 3_1 \# 2^2_1 \rightarrow 3_1 \# 4_1 \rightarrow {4_1^2}'! \rightarrow 0_1 \rightarrow 0^2_1$ is also a parsimonious length 6 pathway.  

\subsubsection{Lower bounds on number of recombinant events:}
Additionally our work can give lower bounds on the minimum number of recombinant events separating DNA configurations.  For example, if both ${4_1^2}'$ and $5^2_1$ arise as products of site-specific recombination, then Corollary \ref{cor:2comp} implies that there must be at least four rounds of recombination between them (as $Lk({4_1^2}'$) = 2, Lk($5^2_1$) = 0, Arf(${4_1^2}'$) = 0 and Arf($5^2_1$) = 1).  This in turn, by Proposition \ref{prop:comp} implies there must exist a 2-component link product $L$ such that $n = d({4_1^2}',L) + d(L,5^2_1)$ and $2 \leq d({4_1^2}', L) \leq n-2$.  

%

\subsubsection{Determining chirality:}

The results above can also determine the chirality of the resulting product knots.  For example the large serine recombinase $\phi$C31 can perform processive recombination.  The resulting products include $2^2_1, 4_1, 5^2_1$ and $6_2$, with the chirality of the last two products unknown. 
From Table \ref{Table:KL} we can see that the chirality of the Whitehead link $5^2_1$ must be the same as the $6_2$ knot, since $d(5^2_1,6_2) (= d(5^2_1!,6_2!)) = 1$ but $d({5_1^{2}}!, 6_2) (= d(5^2_1,6_2!))$ = 3.  This conservation of chirality supports the subunit rotation mechanism discussed in \cite{Femi}.  

As another indicative example, consider the processive recombination events mediated by Tn3 resolvase, as discussed above.  Using electron microscopy, the Cozzarelli group determined that  resulting 6-crossing knot was $6_2!$ \cite{Cozz}.  From our Table \ref{Table:KL} one can see that the minimal coherent pathway between $5^2_1$ and $6_2$ is 1, but between ${5_1^{2}}'$ and $6_2$ is 3.  Thus we can infer that the chirality of the resulting Whitehead links must be ${5^2_1}!$ and not $5^2_1$.

Similarly, when $\lambda$ Int performs recombination on a $2^2_1$ substrate with either so-called PB sites, yielding $4_ 1, 6_1, 8_1$ and $10_1$, or with the LR sites, yielding $0_1, 4_1, 6_1, 8_1, 10_1, 12_1$, then if the chirality of the substrate is known, Table \ref{Table:KL} will determine the chirality of the products.  

\subsubsection{Determining orientation:}

In addition to showing existence of intermediates and/or pathways, our results in Section \ref{s:Tables} can also determine the orientation (i.e. the order of the basepair sequence) of the underlying DNA molecules.  For example, the site-specific recombinases Cre \cite{Hoess}, Flp \cite{Cox} and $\lambda$ Int \cite{Speng} will recombine unknotted circular DNA molecules with direct sites to yield a spectrum of $(2,n)$-torus link products.  Our work in Table \ref{Table:KL} shows that these torus links must be antiparallel, i.e. must be ${(2m)^2_1}'$ and not $(2m)^2_1$  for $m = \{1,2,..,6\}$ for Cre and Flp, and $m = \{4,6\}$ for $\lambda$ Int.


\section*{Appendix}

\end{document}

%% file: DistanceKK1.tex


\begin{tabular}{|c||ccccccccc|}
\hline 
&$3_1$&$4_1$&$5_1$&$5_2$&$6_1$&$6_2$&$6_3$&$3_1\sharp 3_1$&$3_1\sharp 3_1!$\\
\hline
\hline
$0_1$
&$2$($2^2_1$)
&$2$($2^2_1$)
&$4$(I$3_1$)
&$2$($2^2_1$)
&$2$(${2^2_1}'$)
&$2$($2^2_1$)
&$2$($2^2_1$)
&$4$(I$3_1$)
&$2$($3_1\sharp {2^{2}_1}'$)
\\
$3_1$
&$0$
&$2$($2^2_1$)
&$2$($3_1\sharp 2^2_1$)
&$2$($2^2_1$)
&
$4^\dagger$($0_1$)
&$2$($2^2_1$)
&$2$($2^2_1$)
&$2$($3_1\sharp 2^2_1$)
&$2$($3_1\sharp {2^2_1}'$)
\\
$3_1!$
&$4$(I$0_1$)
&$2$(${2^{2}_1}'$)
&$6$(I$0_1$)
&$4$(I$0_1$)
&$2$(${2^{2}_1}'$)
&$4$(I$0_1$)
&$2$(${2^2_1}'$)
&$6$(I$0_1$)
&$2$($3_1!\sharp 2^2_1$)
\\
\hline
$4_1$
&
&$0$
&$4$(I$3_1$)
&$2$($2^2_1$)
&$2$(${2^2_1}'$)
&$2$($2^2_1$)
&$2$($2^2_1$)
&$4$(I$3_1$)
&$4$(II$0_1$)
\\
$5_1$
&
&
&$0$
&$2$($3_1\sharp 2^2_1$)
&$4$(I$3_1$)
&$2$($4^2_1$)
&$4$(I$3_1$)
&$2$($6^2_1$)
&$4$(I$3_1$)
\\
$5_1!$
&
&
&$8$(I$0_1$)
&$6$(I$0_1$)
&$4$(I$3_1!$)
&$6$(I$0_1$)
&$4$(I$3_1!$)
&$8$(I$0_1$)
&$4$(I$3_1!$)
\\
\hline
$5_2$
&
&
&
&$0$
&$2$(${4^2_1}'!$)
&$2$($2^2_1$)
&$2$($2^2_1$)
&$2$($3_1\sharp 2^2_1$)
&$4$(II$0_1$)
\\
$5_2!$
&
&
&
&$4$(I$0_1$)
&$2$(${2^2_1}'$)
&$4$(I$0_1$)
&$2$(${2^2_1}'$)
&$6$(I$0_1$)
&$4$(II$0_1$)
\\
\hline
$6_1$
&
&
&
&
&$0$
&$2$($4_1\sharp 2^2_1$)
&$2$(${2^2_1}'$)
&$4$(I$3_1$)
&$2$($0^2_1$)
\\
$6_1!$
&
&
&
&
&$2$($0^2_1$)
&$2$($2^2_1$)
&$2$($2^2_1$)
&$4$(I$3_1$)
&$2$($0^2_1$)
\\
\hline
$6_2$
&
&
&
&
&
&$0$
&$2$($2^2_1$)
&$4$(II$3_1$)
&$2$($3_1\sharp {2^2_1}'$)
\\
$6_2!$
&
&
&
&
&
&$4$(I$0_1$)
&$2$(${2^2_1}'$)
&$6$(I$0_1$)
&$2$($3_1!\sharp 2^2_1$)
\\
$6_3$
&
&
&
&
&
&
&$0$
&$4$(I$3_1$)
&$2$($3_1\sharp {2^2_1}'$)
\\
\hline
$3_1\sharp 3_1$
&
&
&
&
&
&
&
&$0$
&$4$(I$3_1$)
\\
$3_1!\sharp 3_1!$
&
&
&
&
&
&
&
&$8$(I$0_1$)
&$4$(I$3_1!$)
\\
\hline
\end{tabular}



%% file: DistanceKK2.tex


\begin{tabular}{|c||cccccccc|}
\hline 

&$7_1$&$7_2$&$7_3$&$7_4$&$7_5$&$7_6$&$7_7$&$3_1\sharp 4_1$\\
\hline
\hline
$0_1$
&$6$(I$3_1$)
&$2$($2^2_1$)
&$4$(I$3_1!$)
&$2$(${4^2_1}'$)
&$4$(I$3_1$)
&$2$($2^2_1$)
&$2$($2^2_1$)
&$2$(${4^2_1}'!$)
\\
$3_1$
&$4$(I$5_1$)
&$2$($2^2_1$)
&$6$(I$0_1$)
&$4$(I$0_1$)
&$2$($3_1\sharp 2^2_1$)
&$2$($2^2_1$)
&$2$($2^2_1$)
&$2$($3_1\sharp 2^2_1$)
\\
$3_1!$
&$8$(I$0_1$)
&$4$(I$0_1$)
&$2$($4^2_1!$)
&$2$($3_1!\sharp {2^2_1}'$)
&$6$(I$0_1$)
&$4$(I$0_1$)
&$2$($3_1!\sharp 2^2_1$)
&$4$(I$0_1$)
\\
\hline
$4_1$
&$6$(I$3_1$)
&$2$($2^2_1$)
&$4$(I$3_1!$)
&$2,4$($0_1$)
&$4$(I$3_1$)
&$2$($2^2_1$)
&$2$($2^2_1$)
&$2$($4_1\sharp 2^2_1$)
\\
$5_1$
&$2$($6^2_1$)
&$2$($4^2_1$)
&$8$(I$0_1$)
&$6$(I$0_1$)
&$2$($3_1\sharp 2^2_1$)
&$2$($3_1\sharp 2^2_1$)
&$4$(I$3_1$)
&$2$($3_1\sharp 2^2_1$)
\\
$5_1!$
&$10$(I$0_1$)
&$6$(I$0_1$)
&$2$($4^2_1!$)
&$2$($3_1!\sharp {2^2_1}'$)
&$8$(I$0_1$)
&$6$(I$0_1$)
&$4$(I$3_1!$)
&$6$(I$0_1$)
\\
\hline
$5_2$
&$4$(I$5_1$)
&$2$($2^2_1$)
&$6$(I$0_1$)
&$4$(I$0_1$)
&$2$($3_1\sharp 2^2_1$)
&$2$($3_1\sharp 2^2_1$)
&$2$($2^2_1$)
&$2$(${4^2_1}'!$)
\\
$5_2!$
&$8$(I$0_1$)
&$4$(I$0_1$)
&$2$(${6^2_2}'$)
&$2$(${4^2_1}'$)
&$6$(I$0_1$)
&$4$(I$0_1$)
&$2,4$($0_1$)
&$4$(I$0_1$)
\\
\hline
$6_1$
&$6$(I$3_1$)
&$2,4$($0_1$)
&$4$(I$5_2!$)
&
$4^{\dagger\dagger}(0_1)$
&$4$(I$3_1$)
&$2$($4^2_1$)
&
$4^{\dagger\dagger}(0_1)$
&$2$(${4^2_1}'!$)
\\
$6_1!$
&$6$(I$3_1$)
&$2$($2^2_1$)
&$4$(I$3_1!$)
&$2$(${4^2_1}'$)
&$4$(I$3_1$)
&$2$($2^2_1$)
&$2$($2^2_1$)
&
$4^{\dagger\dagger}(0_1)$
\\
\hline
$6_2$
&$4$(I$5_1$)
&$2$($2^2_1$)
&$6$(I$0_1$)
&$4$(I$0_1$)
&$2,4$($3_1$)
&$2$($2^2_1$)
&$2$($2^2_1$)
&
$2$($3_1\sharp {2^2_1}'$)
\\
$6_2!$
&$8$(I$0_1$)
&$4$(I$0_1$)
&$2$($4^2_1!$)
&$2,4$($0_1$)
&$6$(I$0_1$)
&$4$(I$0_1$)
&$2$($3_1!\sharp 2^2_1$)
&$4$(I$0_1$)
\\
$6_3$
&$6$(I$3_1$)
&$2$($2^2_1$)
&$4$(I$3_1!$)
&$2,4$($0_1$)
&$4$(I$3_1$)
&$2$($2^2_1$)
&$2$($2^2_1$)
&$2$($3_1\sharp {2^2_1}'$)
\\
\hline
$3_1\sharp 3_1$
&$2$($6^2_1$)
&$4$(II$3_1$)
&$8$(I$0_1$)
&$6$(I$0_1$)
&$2$($3_1\sharp 2^2_1$)
&$2$($3_1\sharp 2^2_1$)
&$4$(I$3_1$)
&$2$($3_1\sharp 2^2_1$)
\\
$3_1!\sharp 3_1!$
&$10$(I$0_1$)
&$6$(I$0_1$)
&$4$(II$3_1!$)
&$2$($3_1!\sharp {2^2_1}'$)
&$8$(I$0_1$)
&$6$(I$0_1$)
&$4$(I$0_1$)
&$6$(I$0_1$)
\\
$3_1\sharp 3_1!$
&$6$(I$3_1$)
&$2,4$($0_1$)
&$4$(I$3_1!$)
&$2,4$($0_1$)
&$4$(I$3_1$)
&$4$(II$0_1$)
&$2$($3_1!\sharp 2^2_1$)
&$2$($3_1\sharp {2^2_1}'$)
\\
\hline
$7_1$
&$0$
&$4$(I$5_1$)
&$10$(I$0_1$)
&$8$(I$0_1$)
&$2$($7_1\sharp {2^2_1}'$)
&$4$(I$5_1$)
&$6$(I$3_1$)
&$4$(I$5_1$)
\\
$7_1!$
&$12$(I$0_1$)
&$8$(I$0_1$)
&$2$($7_1!\sharp 2^2_1$)
&$4$(I$5_1!$)
&$10$(I$0_1$)
&$8$(I$0_1$)
&$6$(I$3_1!$)
&$8$(I$0_1$)
\\
\hline
$7_2$
&
&$0$
&$6$(I$0_1$)
&$4$(I$0_1$)
&$2$($7_2\sharp 2^2_1$)
&$2$($2^2_1$)
&$2$($2^2_1$)
&$2,4$($0_1$)
\\
$7_2!$
&
&$4$(I$0_1$)
&$2$($4^2_1!$)
&$2$($9_5\sharp 2^2_1$)
&$6$(I$0_1$)
&$4$(I$0_1$)
&$2,4$($0_1$)
&$4$(I$0_1$)
\\
\hline
$7_3$
&
&
&$0$
&$2$($7_3\sharp 2^2_1$)
&$8$(I$0_1$)
&$6$(I$0_1$)
&$4$(I$3_1!$)
&$6$(I$0_1$)
\\
$7_3!$
&
&
&$8$(I$0_1$)
&$6$(I$0_1$)
&$2$($6^2_2$)
&$2$($4^2_1$)
&$4$(I$3_1$)
&$2,4$($3_1$)
\\
\hline
$7_4$
&
&
&
&$0$
&$6$(I$0_1$)
&$4$(I$0_1$)
&$2,4$($0_1$)
&$4$(I$0_1$)
\\
$7_4!$
&
&
&
&$4$(I$0_1$)
&$2$($3_1\sharp 2^2_1$)
&$2$($3_1\sharp 2^2_1$)
&
$4^{\dagger\dagger}(0_1)$
&$2$(${4^2_1}'!$)
\\
\hline
$7_5$
&
&
&
&
&$0$
&$2$($3_1\sharp 2^2_1$)
&$4$(I$3_1$)
&$2$($3_1\sharp 2^2_1$)
\\
$7_5!$
&
&
&
&
&$8$(I$0_1$)
&$6$(I$0_1$)
&$4$(I$3_1!$)
&$6$(I$0_1$)
\\
\hline
$7_6$
&
&
&
&
&
&$0$
&$2$($2^2_1$)
&$2$($3_1\sharp 2^2_1$)
\\
$7_6!$
&
&
&
&
&
&$4$(I$0_1$)
&$2$($4^2_1\sharp {2^2_1}'$)
&
$4^{\dagger\dagger}(0_1)$
\\
\hline
$7_7$
&
&
&
&
&
&
&$0$
&$2,4$($0_1$)
\\
$7_7!$
&
&
&
&
&
&
&
$4^{\dagger\dagger}(0_1)$
&$2$($3_1\sharp {2^2_1}'$)
\\
$3_1!\sharp 4_1$
&
&
&
&
&
&
&
&$4$(I$0_1$)
\\
\hline
\end{tabular}



%% file: DistanceKL1.tex


\begin{tabular}{|c||cccccccccc|}
\hline 
&$0_1$&$3_1$&$4_1$&$5_1$&$5_2$&$6_1$&$6_2$&$6_3$&$3_1\sharp 3_1$&$3_1\sharp 3_1!$\\
\hline
\hline
$0^2_1$
&$1$
&$3$(I$0_1$)
&$3$(II$0_1$)
&$5$(I$0_1$)
&$3$(I$0_1$)
&$1$
&$3$(I$0_1$)
&$3$(III$0_1$)
&$5$(I$0_1$)
&$1$
\\
$2^2_1$
&$1$
&$1$
&$1$
&$3$(I$3_1$)
&$1$
&$3$(II$0_1$)
&$1$
&$1$
&$3$(I$3_1$)
&$3$(II$0_1$)
\\
${2^2_1}'$
&$1$
&$3$(I$0_1$)
&$1$
&$5$(I$0_1$)
&$3$(I$0_1$)
&$1$
&$3$(I$0_1$)
&$1$
&$5$(I$0_1$)
&$3$(II$0_1$)
\\
\hline
$4^2_1$
&$3$(I$3_1$)
&$1$
&$3$(I$3_1$)
&$1$
&$3$(IV$3_1$)
&$3$(I$3_1$)
&$1$
&$3$(I$3_1$)
&$3$(II$3_1$)
&$3$(I$3_1$)
\\
$4^2_1!$
&$3$(I$3_1!$)
&$5$(I$3_1!$)
&$3$(I$3_1!$)
&$7$(I$3_1!$)
&$5$(I$3_1!$)
&$3$(I$3_1!$)
&$5$(I$3_1!$)
&$3$(I$3_1!$)
&$7$(I$3_1!$)
&$3$(I$3_1!$)
\\
${4^2_1}'$
&$1$
&$3$(I$0_1$)
&$3$(III$0_1$)
&$5$(I$0_1$)
&$3$(I$0_1$)
&$3$(II$0_1$)
&$3$(I$0_1$)
&$3$(IV$0_1$)
&$5$(I$0_1$)
&$3$(II$0_1$)
\\
${4^2_1}'!$
&$1$
&$3$(II$0_1$)
&$3$(III$0_1$)
&$3$(I$5_2$)
&$1$
&$1$
&$3$(IV$0_1$)
&$3$(IV$0_1$)
&$3$(I$5_2$)
&$3$(II$0_1$)
\\
\hline
$5^2_1$
&$3$(IV$3_1$)
&$1$
&$1$
&$3$(I$3_1$)
&$3$(IV$3_1$)
&$3$(II$3_1$)
&$1$
&$1$
&$3$(I$3_1$)
&$3$(II$3_1$)
\\
$5^2_1!$
&$3$(IV$3_1$)
&$3$(I$4_1$)
&$1$
&$5$(I$4_1$)
&$3$(I$4_1$)
&$3$(IV$4_1$)
&$3$(I$4_1$)
&$1$
&$5$(I$4_1$)
&$3$(I$3_1$)
\\
\hline
$3_1\sharp 2^2_1$
&$3$(I$3_1$)
&$1$
&$3$(I$3_1$)
&$1$
&$1$
&$3$(I$3_1$)
&$3$(II$3_1$)
&$3$(I$3_1$)
&$1$
&$3$(I$3_1$)
\\
$3_1!\sharp {2^2_1}'$
&$3$(I$3_1!$)
&$5$(I$3_1!$)
&$3$(I$3_1!$)
&$7$(I$3_1!$)
&$5$(I$3_1!$)
&$3$(I$3_1!$)
&$5$(I$3_1!$)
&$3$(I$3_1!$)
&$7$(I$3_1!$)
&$3$(I$3_1!$)
\\
$3_1\sharp {2^2_1}'$
&$1$
&$1$
&$3$(II$0_1$)
&$3$(I$3_1$)
&$3$(II$0_1$)
&
$3^{\dagger\dagger}(0_1)$
&$1$
&$1$
&$3$(I$0_1$)
&$1$
\\
$3_1!\sharp 2^2_1$
&$1$
&$3$(I$0_1$)
&$3$(II$0_1$)
&$5$(I$0_1$)
&$3$(I$0_1$)
&
$3^{\dagger\dagger}(0_1)$
&$3$(I$0_1$)
&$1$
&$5$(I$0_1$)
&$1$
\\
\hline
$6^2_1$
&$5$(I$5_1$)
&$3$(I$5_1$)
&$5$(I$5_1$)
&$1$
&$3$(I$5_1$)
&$5$(I$5_1$)
&$3$(I$5_1$)
&$5$(I$5_1$)
&$1$
&$5$(I$5_1$)
\\
$6^2_1!$
&$5$(I$5_1!$)
&$7$(I$5_1!$)
&$5$(I$5_1!$)
&$9$(I$5_1!$)
&$7$(I$5_1!$)
&$5$(I$5_1!$)
&$7$(I$5_1!$)
&$5$(I$5_1!$)
&$9$(I$5_1!$)
&$5$(I$5_1!$)
\\
${6^2_1}'$
&$1$
&$5$(I$0_1$)
&$3$(II$0_1$)
&$5$(I$0_1$)
&$3$(I$0_1$)
&
$3^{\dagger\dagger}(0_1)$
&$3$(I$0_1$)
&$1,3$($0_1$)
&$5$(I$0_1$)
&
$3^{\dagger\dagger}(0_1)$
\\
${6^2_1}'!$
&$1$
&$1$
&$3$(II$0_1$)
&$3$(I$0_1$)
&$3$(II$0_1$)
&
$3^{\dagger\dagger}(0_1)$
&
$1,3^{\dagger}(0_1)$
&$1,3$($0_1$)
&$3$(I$3_1$)
&
$3^{\dagger\dagger}(0_1)$
\\
\hline
$6^2_2$
&$3$(I$3_1$)
&$1$
&$3$(I$3_1$)
&$1^{\dagger\dagger\dagger}$
&$1$
&$3$(I$3_1$)
&$3$(III$3_1$)
&$3$(I$3_1$)
&$3$(III$3_1$)
&$3$(I$3_1$)
\\
${6^2_2}'$
&$3$(I$3_1!$)
&$5$(I$3_1!$)
&$3$(I$3_1!$)
&$7$(I$3_1!$)
&$5$(I$3_1!$)
&$3$(I$3_1!$)
&$5$(I$3_1!$)
&$3$(I$3_1!$)
&$7$(I$3_1!$)
&$3$(I$3_1!$)
\\
\hline
$6^2_3$
&$3$(I$5_2$)
&$3$(IV$5_2$)
&$3$(I$5_2$)
&$3$(III$5_2$)
&$1$
&$3$(I$5_2$)
&$3$(II$5_2$)
&$3$(I$5_2$)
&$1$
&$3,5$(I$5_2$)
\\
$6^2_3!$
&$3$(I$5_2!$)
&$5$(I$5_2!$)
&$3$(I$5_2!$)
&$7$(I$5_2!$)
&$5$(I$5_2!$)
&$3$(I$5_2!$)
&$5$(I$5_2!$)
&$3$(I$5_2!$)
&$7$(I$5_2!$)
&$3,5$(I$5_2!$)
\\
${6^2_3}'$
&$3$(I$3_1!$)
&$3$(I$4_1$)
&$1$
&$5$(I$4_1$)
&$3$(I$4_1$)
&$3$(IV$3_1!$)
&$3$(I$4_1$)
&$3$(II$3_1!$)
&$5$(I$4_1$)
&$3$(II$3_1!$)
\\
${6^2_3}'!$
&$3$(II$3_1$)
&$1$
&$1$
&$3$(I$3_1$)
&$3$(IV$3_1$)
&$3$(IV$3_1$)
&$3$(I$3_1$)
&$3$(II$3_1$)
&$3$(I$3_1$)
&$3$(II$3_1$)
\\
\hline
$4_1\sharp 2^2_1$
&$1$
&$3$(II$0_1$)
&$1$
&$3$(I$6_2$)
&$3$(III$0_1$)
&$1$
&$1$
&$3$(III$0_1$)
&$3$(I$7_6$)
&$3$(II$6_2$)
\\
$4_1\sharp {2^2_1}'$
&$1$
&$3$(I$0_1$)
&$1$
&$5$(I$0_1$)
&$3$(I$0_1$)
&$3$(II$0_1$)
&$3$(I$0_1$)
&$3$(III$0_1$)
&$5$(I$0_1$)
&$3$(II$6_2!$)
\\

\hline
\end{tabular}


%% file: DistanceKL2.tex


\begin{tabular}{|c||cccccccc|}
\hline 

&$7_1$&$7_2$&$7_3$&$7_4$&$7_5$&$7_6$&$7_7$&$3_1\sharp 4_1$\\
\hline
\hline
$0^2_1$
&$7$(I$0_1$)
&$3$(I$0_1$)
&$5$(I$0_1$)
&$3$(I$0_1$)
&$5$(I$0_1$)
&$3$(I$0_1$)
&$3$(IV$0_1$)
&$3$(I$0_1$)
\\
$2^2_1$
&$5$(I$3_1$)
&$1$
&$5$(I$0_1$)
&$3$(I$0_1$)
&$3$(I$3_1$)
&$1$
&$1$
&$3$(II$0_1$)
\\
${2^2_1}'$
&$7$(I$0_1$)
&$3$(I$0_1$)
&$3$(I$3_1!$)
&$3$(II$0_1$)
&$5$(I$0_1$)
&$3$(I$0_1$)
&$3$(II$0_1$)
&$3$(I$0_1$)
\\
\hline
$4^2_1$
&$3$(I$5_1$)
&$1$
&$7$(I$3_1$)
&$5$(I$3_1$)
&$3$(IV$3_1$)
&$1$
&$3$(I$3_1$)
&$3$(II$3_1$)
\\
$4^2_1!$
&$9$(I$3_1!$)
&$5$(I$3_1!$)
&$1$
&$3$(II$3_1!$)
&$7$(I$3_1!$)
&$5$(I$3_1!$)
&$3$(I$3_1!$)
&$5$(I$3_1!$)
\\
${4^2_1}'$
&$7$(I$0_1$)
&$3$(I$0_1$)
&$3$(I$5_2!$)
&$1$
&$5$(I$0_1$)
&$3$(I$0_1$)
&$3$(IV$0_1$)
&$3$(I$0_1$)
\\
${4^2_1}'!$
&$5$(I$5_2$)
&$3$(IV$0_1$)
&$5$(I$0_1$)
&$3$(I$0_1$)
&$3$(I$5_2$)
&$3$(I$0_1$)
&$3$(IV$0_1$)
&$1$
\\
\hline
$5^2_1$
&$5$(I$3_1$)
&$1,3$($3_1$)
&$5$(I$4_1$)
&$3,5$(I$3_1$)
&$5$(I$3_1$)
&$1$
&$1$
&$3$(II$3_1$)
\\
$5^2_1!$
&$7$(I$4_1$)
&$3$(I$4_1$)
&$3$(I$3_1!$)
&$3$(II$3_1!$)
&$5$(I$4_1$)
&$3$(I$4_1$)
&$3$(II$3_1!$)
&$3$(IV$4_1$)
\\
\hline
$3_1\sharp 2^2_1$
&$3$(I$5_1$)
&$3$(II$3_1$)
&$7$(I$3_1$)
&$5$(I$3_1$)
&$1$
&$1$
&$3$(I$3_1$)
&$1$
\\
$3_1!\sharp {2^2_1}'$
&$9$(I$3_1!$)
&$5$(I$3_1!$)
&$3$(II$3_1!$)
&$1$
&$7$(I$3_1!$)
&$5$(I$3_1!$)
&$3$(I$3_1!$)
&$5$(I$3_1!$)
\\
$3_1\sharp {2^2_1}'$
&$5$(I$3_1$)
&$1,3$($0_1$)
&$5$(I$0_1$)
&$3$(I$0_1$)
&$3$(I$3_1$)
&$3$(II$0_1$)
&$1,3$($0_1$)
&$1$
\\
$3_1!\sharp 2^2_1$
&$7$(I$0_1$)
&$3$(I$0_1$)
&$3$(I$3_1!$)
&$1,3$($0_1$)
&$5$(I$0_1$)
&$3$(I$0_1$)
&$1$
&$3$(I$0_1$)
\\
\hline
$6^2_1$
&$1$
&$5$(I$5_1$)
&$9$(I$5_1$)
&$7$(I$5_1$)
&$1$
&$3$(I$5_1$)
&$5$(I$5_1$)
&$3$(I$5_1$)
\\
$6^2_1!$
&$11$(I$5_1!$)
&$7$(I$5_1!$)
&$3$(I$5_1!$)
&$3$(I$5_1!$)
&$9$(I$5_1!$)
&$7$(I$5_1!$)
&$5$(I$5_1!$)
&$7$(I$5_1!$)
\\
${6^2_1}'$
&$7$(I$0_1$)
&$3$(I$0_1$)
&$3$(I$3_1!$)
&$1$
&$5$(I$0_1$)
&$3$(I$0_1$)
&$1,3$($0_1$)
&$3$(I$0_1$)
\\
${6^2_1}'!$
&$5$(I$3_1$)
&$1$
&$5$(I$0_1$)
&$3$(I$0_1$)
&$3$(I$3_1$)
&$3$(II$0_1$)
&$1,3$($0_1$)
&$3$(V$0_1$)
\\
\hline
$6^2_2$
&$3$(I$7_3!$)
&$3$(III$3_1$)
&$7$(I$3_1$)
&$5$(I$3_1$)
&$1$
&$3$(III$3_1$)
&$3$(I$3_1$)
&$3$(II$3_1$)
\\
${6^2_2}'$
&$3$(I$3_1!$)
&$5$(I$3_1!$)
&$1$
&$3$(II$3_1!$)
&$7$(I$3_1!$)
&$5$(I$3_1!$)
&$3$(I$3_1!$)
&$5$(I$3_1!$)
\\
\hline
$6^2_3$
&$3$(I$3_1\sharp 3_1$)
&$3$(II$5_2$)
&$7$(I$5_2$)
&$5$(I$5_2$)
&$1$
&$3$(IV$5_2$)
&$3$(I$5_2$)
&$3$(III$5_2$)
\\
$6^2_3!$
&$9$(I$5_2!$)
&$5$(I$5_2!$)
&$3$(II$5_2!$)
&$5$(III$5_2!$)
&$7$(I$5_2!$)
&$5$(I$5_2!$)
&$3,5$(I$5_2!$)
&$5$(I$5_2!$)
\\
${6^2_3}'$
&$7$(I$4_1$)
&$3$(I$4_1$)
&$3$(I$3_1!$)
&$3$(III$3_1!$)
&$5$(I$4_1$)
&$3$(I$4_1$)
&$1$
&$3$(I$4_1$)
\\
${6^2_3}'!$
&$5$(I$3_1$)
&$3$(II$3_1$)
&$5$(I$4_1$)
&$3,5$(I$3_1$)
&$3$(I$3_1$)
&$1$
&$1,3$($3_1$)
&$3$(III$3_1$)
\\
\hline
$4_1\sharp 2^2_1$
&$5$(I$6_2$)
&$1,3$($0_1$)
&$5$(I$0_1$)
&$3$(I$0_1$)
&$3$(I$7_6$)
&$1$
&$3$(II$0_1$)
&$1$
\\
$4_1\sharp {2^2_1}'$
&$7$(I$0_1$)
&$3$(I$0_1$)
&$3$(I$6_2!$)
&$1,3$($0_1$)
&$5$(I$0_1$)
&$3$(I$0_1$)
&$1$
&$3$(I$0_1$)
\\
\hline
\end{tabular}



%% file: DistanceLL1.tex


\begin{tabular}{|c||cccccc|}
\hline 
&$2^2_1$&$4^2_1$&${4^2_1}'$&$5^2_1$&$3_1\sharp 2^2_1$&$3_1\sharp {2^2_1}'$\\
\hline
\hline
$0^2_1$
&$2$($0_1$)
&$4$(I$2^2_1$)
&$2$($0_1$)
&$2$($0_1\sqcup 2^2_1$)*
&$4$(I$2^2_1$)
&$2$($0_1$)
\\
$2^2_1$
&$0$
&$2$($3_1$)
&$2$($0_1$)
&$2$($3_1$)
&$2$($3_1$)
&$2$($0_1$)
\\
${2^2_1}'$
&$2$($0_1$)
&$4$(I$2^2_1$)
&$2$($0_1$)
&$2$($4_1$)
&$4$(I$2^2_1$)
&$2$($0_1$)
\\
\hline
$4^2_1$
&
&$0$
&$4$(I$2^2_1$)
&$2$($3_1$)
&$2$($3_1$)
&$2$($3_1$)
\\
$4^2_1!$
&
&$6$(I$2^2_1$)
&$2$(${2^2_1}'\sharp {2^2_1}'$)*
&$4$(I${2^2_1}'$)
&$6$(I$2^2_1$)
&$4$(I${2^2_1}'$)
\\
${4^2_1}'$
&
&
&$0$
&$4$(IV$2^2_1$)
&$4$(I$2^2_1$)
&$2$($0_1$)
\\
${4^2_1}'!$
&
&
&$2$($0_1$)
&$4$(IV$2^2_1$)
&$2$($5_2$)
&$2$($0_1$)
\\
\hline
$5^2_1$
&
&
&
&$0$
&$2$($3_1$)
&$2$($3_1$)
\\
$5^2_1!$
&
&
&
&$2$($4_1$)
&$4$(I$2^2_1$)
&$2$($6_3$)
\\
\hline
$3_1\sharp 2^2_1$
&
&
&
&
&$0$
&$2$($3_1$)
\\
$3_1!\sharp {2^2_1}'$
&
&
&
&
&$6$(I$2^2_1$)
&$4$(I${2^2_1}'$)
\\
%
$3_1!\sharp 2^2_1$
&
&
&
&
&
&$2$($0_1$)
\\
\hline
\end{tabular}



%% file: DistanceLL2.tex


\begin{tabular}{|c||cccccc|}
\hline 
&$6^2_1$&${6^2_1}'$&$6^2_2$&$6^2_3$&${6^2_3}'$&$4_1\sharp 2^2_1$\\
\hline
\hline
$0^2_1$
&$6$(I$2^2_1$)
&$2$($0_1$)
&$4$(I$2^2_1$)
&$4$(I$2^2_1$)
&$4$(IV$2^2_1$)
&$2$($0_1$)
\\
$2^2_1$
&$4$(I$4^2_1$)
&$2$($0_1$)
&$2$($3_1$)
&$2$($5_2$)
&$2$($4_1$)
&$2$($0_1$)
\\
${2^2_1}'$
&$6$(I$2^2_1$)
&$2$($0_1$)
&$4$(I$2^2_1$)
&$4$(I$2^2_1$)
&$2$($4_1$)
&$2$($0_1$)
\\
\hline
$4^2_1$
&$2$($5_1$)
&$4$(I$2^2_1$)
&$2$($3_1$)
&$2$($2^2_1\sharp 2^2_1$)*
&$4$(I$2^2_1$)
&$2$($7_6$)
\\
$4^2_1!$
&$8$(I$2^2_1$)
&$2$($3_1!$)
&$6$(I$2^2_1$)
&$6$(I$2^2_1$)
&$2$($3_1!$)
&$4$(I$2^2_1$)
\\
${4^2_1}'$
&$6$(I$2^2_1$)
&$2$($0_1$)
&$4$(I$2^2_1$)
&$4$(I$2^2_1$)
&$2$(${4^2_1}'\sharp 2^2_1$)*
&$2$($0_1$)
\\
${4^2_1}'!$
&$4$(I$4^2_1$)
&$2$($0_1$)
&$2$($5_2$)
&$2$($5_2$)
&$4$(IV$2^2_1$)
&$2$($0_1$)
\\
\hline
$5^2_1$
&$6$(I$2^2_1$)
&$2,4$($2^2_1$)
&$2$($3_1!$)
&$4$(IV$2^2_1$)
&$2$($4_1$)
&$2$($4_1$)
\\
$5^2_1!$
&$6$(I$2^2_1$)
&$2$($3_1!$)
&$4$(I$2^2_1$)
&$4$(I$2^2_1$)
&$2$($4_1$)
&$2$($4_1$)
\\
\hline
$3_1\sharp 2^2_1$
&$2$($5_1$)
&$4$(I$2^2_1$)
&$2$($3_1$)
&$2$($5_2$)
&$4$(I$2^2_1$)
&$2$($7_6$)
\\
$3_1!\sharp {2^2_1}'$
&$8$(I$2^2_1$)
&$2$($3_1!$)
&$6$(I$2^2_1$)
&$6$(I$2^2_1$)
&$2$($3_1!$)
&$4$(I${2^2_1}'$)
\\
$3_1\sharp {2^2_1}'$
&$4$(I$4^2_1$)
&$2$($0_1$)
&$2$($3_1$)
&
$4^{\dagger\dagger}(2_1^2)$
&
$4^{\dagger\dagger}(2_1^2)$
&$2$($0_1$)
\\
$3_1!\sharp 2^2_1$
&$6$(I$2^2_1$)
&$2$($0_1$)
&$4$(I$2^2_1$)
&$4$(I$2^2_1$)
&$2$($3_1!$)
&$2$($0_1$)
\\
\hline
$6^2_1$
&$0$
&$6$(I$2^2_1$)
&$2$($4^2_1\sharp 2^2_1$)**
&$2$($3_1\sharp 3_1$)
&$6$(I$2^2_1$)
&$4$(I$4^2_1$)
\\
$6^2_1!$
&$10$(I$2^2_1$)
&$4$(I$4^2_1!$)
&$8$(I$2^2_1$)
&$8$(I$2^2_1$)
&$4$(I$4^2_1!$)
&$6$(I${2^2_1}'$)
\\
${6^2_1}'$
&
&$0$
&$4$(I$2^2_1$)
&$4$(I$2^2_1$)
&$2$($3_1!$)
&$2$($0_1$)
\\
${6^2_1}'!$
&
&$2$($3_1$)
&$2$($3_1$)
&
$4^{\dagger\dagger}(2_1^2)$
&
$4^{\dagger\dagger}(2_1^2)$
&$2$($0_1$)
\\
\hline
$6^2_2$
&
&
&$0$
&$2$($5_2$)
&$4$(I$2^2_1$)
&$2,4$($2^2_1$)
\\
${6^2_2}'$
&
&
&$6$(I$2^2_1$)
&$6$(I$2^2_1$)
&$2$($3_1!$)
&$4$(I${2^2_1}'$)
\\
\hline
$6^2_3$
&
&
&
&$0$
&$4$(I$2^2_1$)
&$2,4$($2^2_1$)
\\
$6^2_3!$
&
&
&
&$6$(I$2^2_1$)
&$2$(${2^2_1}'\sharp {2^2_1}'$)**
&$4$(I${2^2_1}'$)
\\
${6^2_3}'$
&
&
&
&
&$0$
&$2$($4_1$)
\\
${6^2_3}'!$
&
&
&
&
&$2$($4_1$)
&$2$($4_1$)
\\
$4_1\sharp {2^2_1}'$
&
&
&
&
&
&$2$($0_1$)
\\
\hline
\end{tabular}

